\documentclass[a4paper,11pt]{article}
\usepackage[T1]{fontenc}
\usepackage{amsmath, amsfonts, amsthm}
\usepackage{mathrsfs}
\usepackage{xcolor}
\usepackage{bbm}
\usepackage{graphicx}
\usepackage[a4paper]{geometry}
\usepackage{dsfont}
\usepackage{hyperref}
\hypersetup{
colorlinks=true,
linkcolor=blue,
filecolor=magenta,      
urlcolor=cyan
}

\newtheorem{theorem}{Theorem}[section]

\theoremstyle{definition}

\newcommand{\R}{\mathbb{R}}

\DeclareMathOperator{\diag}{Diag}

\newcommand{\ra}{\rightarrow}

\newcommand{\innero}[1]{\langle #1 \rangle}

\newcommand{\ddd}{,\dots,}

\title{{\textbf{ A short proof of Lagrange-Good formula \\ using Dirac delta function }}}
\author{Minh-Toan Nguyen \\ GIPSA-lab, Grenoble Alpes University}
\date{}

\begin{document}
\maketitle
\begin{abstract}
We give a half-page proof of the Lagrange-Good formula, using the Fourier representation of Dirac delta function.
\end{abstract}

\section{Introduction}

Lagrange-Good formula is the multivariate version of the classic Lagrange inversion theorem widely used in enumerative combinatorics. The formula has various proofs using complex analysis \cite{good1960generalizations}, umbral calculus \cite{hofbauer1979short}, operator methods \cite{krattenthaler1988operator}, quantum field theory \cite{abdesselam2002physicist} and symbolic integral \cite{abdesselam2002feynman}. There are also proofs based on combinatorial arguments \cite{gessel1987combinatorial} and  on a determinant formula \cite{huang2017determinant}. We give here a short proof of the formula using Dirac delta function.

The formula, originally proved for complex analytic functions in \cite{good1960generalizations}, holds more generally for formal power series \cite{tutte1975elementary}. This version of the formula is stated as follows.

\begin{theorem}
Consider formal power series of $n$ variables $x_1 \ddd x_n$, with coefficients in $\R$ or $\mathbb C$. Denote $x=(x_1 \ddd x_n)$.
Consider $n$ formal power series $f_i(x)$ for $i =1\ddd n$. Let $g_i(x)$ be the formal power series that solve the following system of equations
\begin{align}
	g_i = x_i f_i(g_1 \ddd g_n), \quad i = 1 \ddd n.
\end{align}
Then for any formal power series $\phi$, the coefficient of $[x_1^{k_1} \dots x_n^{k_n}]$ in
\begin{align}\label{lhs}
	\frac{\phi( g( x))}{\det (\delta_{ij} - x_i \partial_j f_i(g) )}
\end{align}
is equal to the coefficient of $[x_1^{k_1} \dots x_n^{k_n}]$ in
\begin{align*}
	\phi(x) f_1(x)^{k_1} \dots f_n(x)^{k_n}.
\end{align*}
\end{theorem}

To prove results about formal power series, we will work with functions that are smooth, analytic in a neighborhood of zero and have a compact support. We denote this class of functions by $\mathcal C$. Algebraic formulas concerning the Taylor coefficients of these function around zero are also valid for formal power series. This is because the regularity imposed on functions in $\mathcal C$ has no impact on the algebraic aspect of the result.

\textbf{Notation.} $\diag(x)$ is the diagonal matrix with the diagonal $x$. $[m]f$ is the coefficient of the monomial $m$ in the series $f$. $C_c^\infty(\R^n)$ is the class of smooth functions from $\R^n$ to $\R$ with compact support.

\section{Dirac delta function}
Let us now recall some basic facts about the Dirac delta function, which is defined as the generalized function \cite{strichartz2003guide} such that
\begin{align}
	\int \delta(x) \phi(x) dx = \phi(0)
\end{align}
for all $\phi \in C_c^\infty(\R^n)$. First, we have
\begin{align}\label{ok}
	\delta(f(x)) = \sum_{ x_i: f(x_i)=0} \frac{ \delta(x-x_i) }{|\det J_{f}(x_i)|}
\end{align}
where $J_{f}(x)$ is the Jacobian matrix of $f$ at $x$,
\begin{align*}
	J_{f}( x) =  [\partial_j  f_i(x)]_{i,j=1}^n
\end{align*}

The Dirac delta function has the following Fourier representation
\begin{align}
	\delta( u) &= \frac{1}{(2\pi)^n}\int_{\R^n} e^{-i \innero{ \xi,  u}} d \xi
\end{align}
which can be written compactly as
\begin{align}\label{drep}
	\delta( u) = \int \hat d v \, e^{-\innero{v,  u}}
\end{align}
where
\begin{align*}
\int \hat d v = \frac{1}{(2\pi i)^n} \int_{i \R^n} dv
\end{align*}
For any function $h \in C_c^\infty(\R^n)$, we have
\begin{align*}
h(x) = \int du \, \delta(u-x)h(u)
\end{align*}
By (\ref{drep}), we have
\begin{align*}
h(x) = \int \hat dv du  \, e^{-\innero{v,u-x}} h(u)
\end{align*}
Applying on both sides the operator $g(\partial)$ where $g$ is some polynomial, we obtain the formula
\begin{align}\label{pk}
	g( \partial)  h(x)|_{x=0} = \int \hat d v du  \, e^{-\innero{v,u}}   g(v) h(u) 
\end{align}
which relates operator calculus with complex integrals.

\section{The proof} \label{proof}

Let $f$ be a function in $\mathcal C$. If $x \in \R^n$ is small enough such that the mapping $u \ra \diag(x) f(u)$ is a contraction, then the equation $u - \diag(x) f(u)$ has unique solution. So there exists $\epsilon >0$ such that for all $x \in B(0, \epsilon)$, the equation $u - \diag(x) f(u)$ has unique solution, denoted by $g(x)$. Let $\phi \in \mathcal C$ and
\begin{align*}
	I(x) = \int d u\,  \phi( u) \delta[  u - \diag(x)  f( u) ], \quad x \in B(0, \epsilon).
\end{align*}
From (\ref{ok}), we have
\begin{align*}
I(x) = \frac{\phi( g( x))}{\det (\delta_{ij} - x_i \partial_j f_i(g) )}.
\end{align*}
From the Fourier representation of Dirac delta function, we have
\begin{align*}
	 I( x) =  \int \hat dv du   \, \phi(u) e^{ -\innero{v,u} + \innero{v, \diag(x) f(u)} }.
\end{align*}
Applying $ \frac{\partial_{x_1}^{k_1} \dots \partial_{x_n}^{k_n}}{k_1! \dots k_n!} $ at $x=0$ on both sides, we obtain
\begin{align*}
	[x_1^{k_1} \dots x_n^{k_n}] I( x) &= \int \hat dv du  \, \phi(u) e^{ -\innero{v, u} }  \frac{(v_1 f_1( u))^{k_1} \dots (v_n f_n( u))^{k_n} }{k_1! \dots k_n!} \\
	&= \frac{\partial_{x_1}^{k_1} \dots \partial_{x_n}^{k_n}}{k_1! \dots k_n!} \phi(x) f_1(x)^{k_1} \dots f_n(x)^{k_n}|_{x=0} \quad \text{(by (\ref{pk}))} \\
	&= [x_1^{k_1} \dots x_n^{k_n}] \phi(x) f_1(x)^{k_1} \dots f_n(x)^{k_n},
\end{align*}
which is the Lagrange-Good formula.

\bibliographystyle{siam}
\bibliography{lagrange-good.bib}

\begin{thebibliography}{1}

\bibitem{abdesselam2002physicist}
{\sc A.~Abdesselam}, {\em {A physicist's proof of the Lagrange-Good
  multivariable inversion formula}}, arXiv preprint math/0208174,  (2002).

\bibitem{abdesselam2002feynman}
\leavevmode\vrule height 2pt depth -1.6pt width 23pt, {\em Feynman diagrams in
  algebraic combinatorics}, arXiv preprint math/0212121,  (2002).

\bibitem{gessel1987combinatorial}
{\sc I.~M. Gessel}, {\em {A combinatorial proof of the multivariable Lagrange
  inversion formula}}, Journal of Combinatorial Theory, Series A, 45 (1987),
  pp.~178--195.

\bibitem{good1960generalizations}
{\sc I.~J. Good}, {\em {Generalizations to several variables of Lagrange's
  expansion, with applications to stochastic processes}}, in Mathematical
  Proceedings of the Cambridge Philosophical Society, vol.~56, Cambridge
  University Press, 1960, pp.~367--380.

\bibitem{hofbauer1979short}
{\sc J.~Hofbauer}, {\em {A short proof of the Lagrange-Good formula}}, Discrete
  Mathematics, 25 (1979), pp.~135--139.

\bibitem{huang2017determinant}
{\sc J.~Huang and X.~Ma}, {\em {A determinant identity implying the
  Lagrange-Good inversion formula}}, Proceedings of the Edinburgh Mathematical
  Society, 60 (2017), pp.~165--176.

\bibitem{krattenthaler1988operator}
{\sc C.~Krattenthaler}, {\em {Operator methods and Lagrange inversion: a
  unified approach to Lagrange formulas}}, Transactions of the American
  Mathematical Society, 305 (1988), pp.~431--465.

\bibitem{strichartz2003guide}
{\sc R.~S. Strichartz}, {\em {A guide to distribution theory and Fourier
  transforms}}, World Scientific Publishing Company, 2003.

\bibitem{tutte1975elementary}
{\sc W.~Tutte}, {\em {On elementary calculus and the Good formula}}, Journal of
  Combinatorial Theory, Series B, 18 (1975), pp.~97--137.

\end{thebibliography}

\end{document}